\documentclass[abstracton,a4paper]{scrartcl}

\usepackage[headsepline]{scrpage2}
\usepackage[latin1]{inputenc}
\usepackage{amsfonts}
\usepackage{amsmath}
\usepackage{amssymb}
\usepackage{amsthm}

\usepackage{url}

\usepackage{graphicx}
\usepackage{pdftricks}
  \begin{psinputs}
    \usepackage{pstricks}
    \usepackage{pst-plot}    
  \end{psinputs} 

\usepackage[pdftex,colorlinks,breaklinks,linkcolor=black,citecolor=black,filecolor=black,menucolor=black,urlcolor=black,pdfauthor={Lizhen Ji and Andreas Weber},pdftitle={Lp spectrum and heat dynamics of locally symmetric spaces with higher rank}, plainpages=false, pdfpagelabels, bookmarksnumbered=true]{hyperref}

\ihead{L. Ji, A. Weber: $L^p$ spectrum and heat dynamics of locally symmetric spaces}

\newtheorem{theorem}{Theorem}[section]
\newtheorem{lemma}[theorem]{Lemma}
\newtheorem{proposition}[theorem]{Proposition}
\newtheorem{corollary}[theorem]{Corollary}
\newtheorem{conjecture}[theorem]{Conjecture}

\theoremstyle{definition}
\newtheorem{definition}[theorem]{Definition}

\theoremstyle{remark}
\newtheorem{remark}[theorem]{\bf Remark}

\newcommand{\Hy}{\mathbb{H}}

\newcommand{\Q}{\mathbb{Q}}
\newcommand{\R}{\mathbb{R}}
\newcommand{\Z}{\mathbb{Z}}
\newcommand{\C}{\mathbb{C}}

\newcommand*\e{\mathrm{e}}

\newcommand*\re{\mathrm{Re}}
\newcommand*\im{\mathrm{Im}}
\newcommand*\id{\mathrm{id}}

\newcommand*\grad{\mathrm{grad}}
\newcommand*\dive{\mathrm{div}}

\newcommand*\isom{\mathrm{Isom}}
\newcommand*\vol{\mathrm{vol\,}}

\newcommand*\rank{\mathrm{rank}}
\newcommand*\qrank{\Q\mbox{-}\mathrm{rank}}

\newcommand*\dom{\mathrm{dom}}
\newcommand*\tr{\mathrm{tr}}
\newcommand*\spa{\mathrm{span}}

\newcommand\ad{\mathrm{ad}}
\newcommand\Ad{\mathrm{Ad}}

\newcommand*\glnc{\mathrm{\it GL}(n,\C)}

\newcommand\bG{{\bf G}}

\newcommand\bP{{\bf P}}
\newcommand\bN{{\bf N}}
\newcommand\bL{{\bf L}}
\newcommand\bS{{\bf S}}
\newcommand\bM{{\bf M}}

\newcommand*\DMp{\Delta_{M,p}}                 
\newcommand*\DMq{\Delta_{M,q}}
\newcommand*\DM{\Delta_M}

\title{$\boldsymbol{L^p}$  spectrum and heat dynamics of locally symmetric spaces of higher rank}
\author{Lizhen Ji\footnote{Email: lji@umich.edu, 
					Address: 1834 East Hall, Ann Arbor, MI 48109-1043, USA.
		     Partially supported by NSF grant DMS 0905283 }\\ 
					{\large Department of Mathematics, University of Michigan}
\and Andreas Weber\footnote{    Email: andreasweber.mail@gmail.com,
						   Address:  Kaiserstr. 89 - 93, 76128 Karlsruhe, Germany.}\\
						   {\large  Institut f\"ur Algebra und Geometrie,
						   KIT} }

\date{}

\begin{document}

\maketitle 
\begin{abstract} 
 The aim of this paper is to study the spectrum of the $L^p$ Laplacian and the dynamics of the $L^p$ heat semigroup on non-compact locally symmetric spaces of higher rank.  Our work here generalizes previously obtained results in the setting
of locally symmetric spaces of rank one to higher rank spaces.
Similarly as in the rank one case, it turns out that the $L^p$ heat semigroup on $M$ has a certain chaotic behavior if $p\in(1,2)$ whereas for $p\geq 2$ such a chaotic behavior never occurs. \\

 \noindent{\em Keywords:} Locally symmetric spaces, $L^p$ heat semigroups, 
          Eisenstein series, $L^p$ spectrum, chaotic semigroups. 
\end{abstract}		

\section{Introduction}

The aim of this paper is to study the spectrum of the $L^p$ Laplacian $\DMp$ and the dynamics of the $L^p$ heat semigroup $e^{-t\DMp}: L^p(M)\to L^p(M)$ on non-compact locally symmetric spaces 
$M=\Gamma\backslash X$ of higher rank with finite volume. 
More precisely, $X$ is a symmetric space of non-comact type and $\Gamma$ a non-uniform arithmetic subgroup of isometries.

The $L^2$ spectral theory is of fundamental importance for locally symmetric spaces
and has been extensively studied in the past (see \cite{MR2189882,0926.11034,MR1025165}
and the references therein).
However, in contrast to the $L^2$ spectrum of the Laplace-Beltrami operator on locally symmetric
spaces, the $L^p$ spectrum, $p\in [1,\infty)$, is only known in special situations.
But there are several reasons to study also the $L^p$ spectrum for $p\neq 2$.\\ 
Firstly, from a physical point of view, the natural space to study heat diffusion is $L^1$:
If the function $u(t,\cdot)\geq 0$ denotes the heat distribution at the time $t$, the total amount of heat 
in some region $\Omega$ is given by the $L^1$ norm of $u(t,\cdot)|_{\Omega}$ and hence, the $L^1$ norm has a physical meaning. But, on the other hand, $L^1$ is more difficult to handle than the reflexive $L^p$ spaces ($p>1$) as the heat semigroup on $L^p$ ($p>1$) is  always bounded analytic whereas this is in general not true for the heat semigroup on $L^1$.\\ 
Secondly, there are already many results for differential operators on domains of 
euclidean space concerning various aspects of $L^p$ spectral theory, see e.g.
\cite{MR1269649,MR1103113,MR1345724,MR1458713,MR836002,MR869525,MR1602267,MR1673414,MR1824257,MR1224619,MR1420468}. 
But, in contrast to the euclidean situation, there are not many examples of manifolds 
whose $L^p$ spectral theory is well understood. This paper continues our work on
$L^p$ spectral theory of locally symmetric spaces in \cite{Ji:yq}, now concentrating on the more
involved higher rank case. To a certain extend, we will reveal the structure of the
$L^p$ spectrum and provide bounds for the $L^p$ spectrum.
A surprising consequence of these results is that 
the $L^p$ heat semigroup has a certain chaotic behavior if $p\in(1,2)$ whereas for $p\geq 2$ such a chaotic behavior never occurs (see Theorem \ref{thm dynamics}). \\

Let us briefly recall related previous work before we give a more detailed description of the contents
of this paper.\\
Davies, Simon, and Taylor completely determined in \cite{MR937635} the $L^p$ spectrum of the Laplacian on the $n$-dimensional hyperbolic space $X=\Hy^n$ and on non-compact quotients 
$M=\Gamma\backslash X$ of $X$ by a geometrically finite group $\Gamma$ of isometries
such that $M$ has no cusps or has finite volume.\\
Taylor generalized some of these results to symmetric spaces $X$ of 
non-compact type  \cite{MR1016445}. More precisely, he proved the following 
result (for a definition of $\rho$ we refer to Section \ref{symmetric spaces}) .
\begin{theorem}[cf. \cite{MR1016445}]\label{taylor}
	Let $X$ denote a symmetric space of non-compact type. Then  for any 
	$p\in [1,\infty)$ we have 
	$$\sigma(\Delta_{X,p}) = {\cal P}_{X,p},$$
	where
	$${\cal P}_{X,p}=\left\{ ||\rho||^2 - z^2 : z\in\C, |\re z|\leq ||\rho||\cdot |\frac{2}{p}-1|\right\}.$$
	Furthermore, if $p>2$ every point in the interior of the parabolic region ${\cal P}_{X,p}$ is an 
	eigenvalue for $\Delta_{X,p}$ and eigenfunctions corresponding to these eigenvalues
    are given by spherical functions.
\end{theorem}
In the case $p\leq 2$ the Helgason-Fourier-Transform ${\cal F}$ turns the Laplacian into a multiplication
operator, i.e. $({\cal F}\Delta f)(\lambda) = 
 	(||\rho||^2 + \langle \lambda,\lambda\rangle){\cal F}f$. 
Hence, together with the $L^p$ inversion formula by Stanton and Tomas \cite{MR518528},
it follows that there are no eigenvalues in the case $p\leq 2$. 	\\
Taylor, using similar ideas as in  \cite{MR937635}, was also able to give an upper bound for 
the $L^p$ spectrum on quotients $M=\Gamma\backslash X$:
\begin{proposition}[Proposition 3.3 in \cite{MR1016445}]\label{taylor upper}
Let $M$ denote a non-compact locally symmetric space whose universal cover is a symmetric
space of non-compact type and assume
$$
 \sigma(\Delta_{M,2}) \subset \{\lambda_0,\ldots,\lambda_r\}\cup [b,\infty),
$$
where $\lambda_j, j=0,\ldots, r$, are eigenvalues of finite multiplicity. Then,
if $\vol(M)<\infty$ or if the injectivity radius of $M$ is bounded away
from $0$, we have for $p\in (1,\infty)$
$$
 \sigma(\Delta_{M,p}) \subset \{\lambda_0,\ldots,\lambda_r\}\cup {\cal P}_{M,p}',
$$
where
$$
 {\cal P}_{M,p}' = \left\{ b - z^2 : z\in\C, |\re z|\leq ||\rho||\cdot |\frac{2}{p}-1|\right\}.
$$
\end{proposition}
Some of these results have been generalized in \cite{Ji:2007fk,Ji:yq,MR2342629} to rank one locally
symmetric spaces. For results concerning the $L^p$ heat dynamics on symmetric spaces
of non-compact type we refer to  \cite{Ji:nr}.\\

In this paper, after recalling some basic definitions and facts concerning $L^p$ heat semigroups
and locally symmetric spaces in the next section, 
we study Eisenstein series $E(\bP|\varphi,\Lambda)$ on 
$M=\Gamma\backslash X$ which are generalized (smooth) eigenfunctions for the 
$L^2$ Laplacian $\Delta_{M,2}$. In order to show that if $p<2$ 
these Eisenstein series are contained in $L^p(M)$  for certain choices of $\Lambda$
(cf. Theorem \ref{thm general parabolic} and Corollary \ref{eigenvalues minimal}), and hence are
eigenfunctions of $\DMp$,
we first derive an  upper bound  for Eisenstein series (Proposition \ref{upper bound}, 
Corollary \ref{corollary eisenstein upper bound minimal}). This bound might be expected by a
reader who is familiar with the theory of Eisenstein series. However, since  we
could not find an upper bound in the form we present it here in the current literature, we included a proof.\\
From these results it then follows that in the case $p\in(1,2)$ there is an open subset of $\C$
consisting of eigenvalues of the $L^p$ Laplacian $\DMp$ whereas in the case $p\geq 2$
the point spectrum of the $L^p$ Laplacian is discrete:\\

{\noindent \bf Theorem \ref{theorem eigenvalues}.}{\it\, 
	Let $\bP$ denote a minimal rational parabolic subgroup of $\bG$ and
	\begin{equation}
 		{\cal P}_{M,p}(\rho_{\bP}) = \left\{ ||\rho_{\bP}||^2 - z^2 : 
 				z\in\C, |\re z|\leq ||\rho_{\bP}||\cdot |\frac{2}{p}-1|\right\}\subset \C.
	\end{equation} 	
	\begin{itemize}
		\item[\textup{(a)}]  If $p\in (1,2)$, there exists a discrete set 
			$B\subset {\cal P}_{M,p}(\rho_{\bP})$ such that all points in the interior of 
			${\cal P}_{M,p}(\rho_{\bP}) \setminus B$ are eigenvalues of $\DMp$.
		\item[\textup{(b)}] If $p\in [2,\infty)$, the point spectrum $\sigma_{pt}(\DMp)$ is a 
			discrete subset of $[0,\infty)$.		
	\end{itemize}
	Furthermore, ${\cal P}_{M,p}(\rho_{\bP}) \subset \sigma(\DMp)$ for all $p\in[1,\infty)$.\\
}

Both of these facts are a main ingredient in the proof of  Theorem \ref{thm dynamics} which is
devoted to the dynamics of the $L^p$ heat semigroup
$e^{-t\DMp}: L^p(M)\to L^p(M):\\$

{\noindent \bf Theorem \ref{thm dynamics}.} {\it\,
	Let $M=\Gamma\backslash X$ denote a non-compact locally symmetric space with
	arithmetic fundamental group $\Gamma$.
 	\begin{itemize}
   		\item[\textup{(a)}]
   			If $p\in (1,2)$ there is a constant $c_p>0$ such that for any $c>c_p$ the semigroup
   			$$
    			e^{-t(\DMp -c)}: L^p(M) \to L^p(M)
   			$$
   			is subspace chaotic.
   		\item[\textup{(b)}]
   			For any $p\geq 2$ and $c\in\R$ the semigroup $e^{-t(\DMp -c)}: L^p(M) \to L^p(M)$ 
   			is not subspace chaotic.
 		\end{itemize}
}
Note that Theorem \ref{theorem eigenvalues} seems to be the first result in which it is shown that some subset of $\C$ with non-empty interior (a parabolic region) is contained in the $L^p$ spectrum of the Laplacian on a locally symmetric space of higher $\Q$-rank.
Nevertheless, the exact determination of the $L^p$ spectrum of a non-compact locally
symmetric space with finite volume remains open in the higher rank case 
(cf. Conjecture \ref{conjecture}).
Even for rank one spaces the $L^p$ spectrum is completely 
known only in certain situations, cf. \cite{MR937635,Ji:yq}.
However, as the $L^p$ spectrum of a Riemannian product $M=M_1\times M_2$
equals the set theoretic sum of the $L^p$ spectra of its factors, i.e.
$\sigma(\DMp) = \sigma(\Delta_{M_1,p}) + \sigma(\Delta_{M_2,p})$, cf. \cite{Weber:2008ve},
we can restric ourselves to irreducible manifolds. 
%
\section{Preliminaries}
%
\subsection{The heat semigroup on $\boldsymbol{L^p}$ spaces}\label{heat semigroup}

In this section we denote by $M$ an arbitrary complete Riemannian manifold and by
$\Delta=-\dive(\grad)$ the Laplace-Beltrami operator acting on differentiable functions of $M$.
If we denote by $\DM$ the Laplacian on $L^2(M)$ with domain $C_c^{\infty}(M)$ (the set of differentiable functions with compact support), this is an essentially self-adjoint operator and hence, 
its closure $\Delta_{M,2}$ is a self-adjoint operator on the Hilbert space $L^2(M)$. 
Since $\Delta_{M,2}$ is positive,
$-\Delta_{M,2}$ generates a bounded analytic semigroup $e^{-t\Delta_{M,2}}$ on $L^2(M)$. 
The semigroup $e^{-t\Delta_{M,2}}$
is a {\em submarkovian semigroup} (i.e., $e^{-t\Delta_{M,2}}$ is positive and a contraction on 
$L^{\infty}(M)$
for any $t\geq 0$) and we therefore have the following:
\begin{itemize}
\item[(1)] The semigroup $e^{-t\Delta_{M,2}}$ leaves the set $L^1(M)\cap L^{\infty}(M)\subset L^2(M)$ 
		invariant and thus,
		$e^{-t\Delta_{M,2}}|_{L^1\cap L^{\infty}}$ may be extended to a positive contraction semigroup
		$T_p(t)$ on $L^p(M)$ for any $p\in [1,\infty]$.
		 These semigroups are strongly continuous if $p\in [1,\infty)$ and {\em consistent}
		 in the sense that $T_p(t)|_{L^p\cap L^q} = T_q(t)|_{L^p\cap L^q}$. 
\item[(2)] Furthermore, if $p\in (1,\infty)$, the semigroup $T_p(t)$ is a bounded analytic semigroup
		with angle of analyticity $\theta_p \geq \frac{\pi}{2} - \arctan\frac{|p-2|}{2\sqrt{p-1}}$.	
\end{itemize} 
For a proof of (1) we refer to \cite[Theorem 1.4.1]{MR1103113}. For (2) see \cite{MR1224619}.
In general, the semigroup $T_1(t)$ needs not be analytic. However, if $M$ has bounded geometry
$T_1(t)$ is analytic in {\em some} sector (cf. \cite{MR924464,MR1023321}).

In the following, we denote by $-\DMp$ the generator of $T_p(t)$ and by $\sigma(\DMp)$ the spectrum of $\DMp$. Furthermore, we write
$e^{-t\DMp}$ for the semigroup $T_p(t)$.
Because of (2) from above, 
the $L^p$ spectrum $\sigma(\DMp)$ of $M$ is contained in the sector 
\begin{multline*}
\left\{ z\in \C\setminus\{0\} : |\arg(z)| \leq \frac{\pi}{2}-\theta_p\right\}\cup\{0\} \subset\\
     \left\{ z\in \C\setminus\{0\} : |\arg(z)| \leq \arctan\frac{|p-2|}{2\sqrt{p-1}} \right\}\cup\{0\} = \Sigma_p.
\end{multline*}     
If we identify as usual the dual space of $L^p(M), 1\leq p<\infty$, with 
$L^{p'}(M), \frac{1}{p}+\frac{1}{p'}=1$, the dual operator of $\DMp$ equals $\Delta_{M,p'}$
and therefore we always have $\sigma(\DMp) = \sigma(\Delta_{M,p'})$. 
It should also be mentioned that the family $\DMp, p\geq 1,$ is consistent, which means that
the restrictions of $\DMp$ and $\Delta_{M,q}$ to the intersection of their domains coincide:
\begin{lemma}\label{lemma consistent}
If $p,q \in [1,\infty)$, the operators $\DMp$ and $\Delta_{M,q}$ are consistent, i.e.
$$ 
\DMp f = \Delta_{M,q} f\qquad\mbox{for any~} f\in \dom(\DMp)\cap\dom(\Delta_{M,q}).
$$
\end{lemma}
For a proof see e.g. \cite[Lemma 2.1]{Ji:yq}.\\

Since it is not obvious that a differentiable $L^p$ function $f$ that 
satisfies the eigenequation $\Delta f = \mu f$ is  contained in the domain of $\DMp$, we
state this result in a lemma:
\begin{lemma}\label{lemma Lp eigenfunctions}
 Let $p\in (1,\infty)$ and $f: M\to \R$ denote a differentiable function such that $f\in L^p(M)$ and 
 $\Delta f = \mu f$ for some $\mu\in\R$.  Then $f\in \dom(\DMp)$ and $\DMp f = \mu f$.
\end{lemma}
For a proof we refer to \cite[Corollary 2.3]{Ji:yq}.\\
An essential ingredient in this proof is  a uniqueness result concerning $L^p$ solutions of the heat equation by Strichartz, cf. \cite[Theorem 3.9]{MR705991}.
As in this theorem $p$ is required to be contained in $(1,\infty)$, we have excluded the case $p=1$ here as well. 
%
\subsubsection{Manifolds with finite volume}

If $M$ is Riemannian manifold with finite volume we have by H\"older's inequality
$L^q(M) \hookrightarrow L^p(M)$ for $1\leq p\leq q\leq \infty$. Hence, by consistency, the 
semigroup $e^{-t\Delta_{M,p}}$ can be regarded as extension of the semigroup
$e^{-t\Delta_{M,q}}$, $p\leq q$. It also follows an analogous result
for the $L^p$ Laplacians, which is stronger than Lemma \ref{lemma consistent} in the case of 
finite volume:
\begin{lemma}\label{domains}
 Let $M$ denote a complete Riemannian manifold with finite volume. If $1\leq p\leq q<\infty$, we have
 $$
  \dom(\Delta_{M,q})\subset\dom(\DMp)
 $$
 and for $f\in\dom(\Delta_{M,q})$ it follows $\Delta_{M,q}f=\DMp f$, i.e. $\DMp$ is an extension of
 $\Delta_{M,q}$.
\end{lemma}
\begin{proof}
Let $f\in\dom(\Delta_{M,q})$. Because of H\"older's inequality and because of consistency of
the $L^p$ heat semigroups, we have
 \begin{eqnarray*}
  || \frac{1}{t}(e^{-\DMp}f - f) - \Delta_{M,q}f ||_{L^p}  
  	& \leq & C  || \frac{1}{t}(e^{-\DMp}f - f) - \Delta_{M,q}f ||_{L^q} \\
	& = & C  || \frac{1}{t}(e^{-\Delta_{M,q}}f - f) - \Delta_{M,q}f ||_{L^q}\\
	&\to& 0 \qquad (t\to 0),
 \end{eqnarray*}
i.e. $f\in \dom(\DMp)$ with $\Delta_{M,q}f=\DMp f$.
\end{proof}
\begin{proposition}
 Let $M$ denote a complete Riemannian manifold with finite volume. 
 If $2\leq p\leq q<\infty$, we have
 $$
  \sigma(\DMp) \subset \sigma(\Delta_{M,q}).
 $$
 In particular, for $p\in (1,\infty)$,
 it follows that
 $$
 \sigma(\Delta_{M,2}) \subset \sigma(\DMp).
 $$
\end{proposition}
\begin{proof}
For a proof of the first statement we refer to \cite[Proposition 3.3]{MR2342629}.
The second statement follows from the first together with the fact
$\sigma(\DMp)=\sigma(\Delta_{M,p'}), \frac{1}{p} + \frac{1}{p'}=1$.
\end{proof}
Note that in the case of a compact manifold $M$, we always have
$\sigma(\DMp) = \sigma(\Delta_{M,2})$ whereas for non-compact manifolds this needs
not be true \cite{MR1250269}. Examples for manifolds whose spectrum depends
non-trivially on $p$ are non-compact arithmetic locally symmetric spaces, cf. Section \ref{Lp spectrum}.
It is actually shown there that the $L^p$ spectrum (if $p\neq 2$) contains a subset of $\C$
with non-empty interior (see also \cite{MR937635,Ji:nr,Ji:yq}).

On the one hand, this proposition states that the $L^p$ spectrum is bigger than the $L^2$
spectrum. On the other hand, how much bigger it is, is often difficult to say. In Conjecture \ref{conjecture}
we give a precise picture of the $L^p$ spectrum of non-compact locally symmetric spaces $M=\Gamma\backslash X$ with arithmetic $\Gamma$.

\subsection{Locally symmetric spaces}\label{symmetric spaces}
Let $X$ denote always a symmetric space of non-compact type. Then
$G= \isom^0(X)$ is a non-compact, semi-simple Lie group with trivial center 
that acts transitively on $X$ and $X\cong G/K$ for any maximal compact subgroup 
$K$ of $G$. We denote
the respective Lie algebras by $\mathfrak{g}$ and $\mathfrak{k}$. Given a corresponding Cartan
involution $\theta: \mathfrak{g}\to\mathfrak{g}$ we obtain the Cartan decomposition
$\mathfrak{g}=\mathfrak{k}\oplus\mathfrak{p}$ of $\mathfrak{g}$ into the eigenspaces of $\theta$. 
The subspace
$\mathfrak{p}$ of $\mathfrak{g}$ can be identified with the tangent space $T_{eK}X$ and we assume
that the Riemannian metric $\langle\cdot,\cdot\rangle$ of $X$ in $\mathfrak{p}\cong T_{eK}X$ 
coincides with the restriction of the Killing form 
$B(Y,Z) = \tr(\ad Y\circ \ad Z ), Y, Z\in \mathfrak{g},$ to $\mathfrak{p}$. 

For any maximal abelian subspace $\mathfrak{a}\subset \mathfrak{p}$ we refer to 
$\Phi(\mathfrak{g},\mathfrak{a})$ as the set of restricted roots for the pair $(\mathfrak{g},\mathfrak{a})$,
i.e. $\Phi(\mathfrak{g},\mathfrak{a})$  contains all $\alpha\in \mathfrak{a}^*\setminus\{0\}$ such that
$$ 
\mathfrak{h}_{\alpha} = 
	\{ Y\in \mathfrak{g} : \ad(H)(Y) = \alpha(H)Y \mbox{~for all~} H\in\mathfrak{a} \}\neq \{0\}.
$$
These subspaces $ \mathfrak{h}_{\alpha}\neq \{0\}$ are called root spaces.\\
Once a positive Weyl chamber $\mathfrak{a}^+$ in $\mathfrak{a}$ is chosen, we denote by
$\Phi(\mathfrak{g},\mathfrak{a}) ^+$ the  subset of positive roots and by 
$$
	\rho = \frac{1}{2}\sum_{\alpha\in\Phi^+(\mathfrak{g},\mathfrak{a}) } 
	(\dim \mathfrak{h}_{\alpha})\alpha
$$ 
half the sum of the positive roots (counted according to their multiplicity).\\

In what follows, we denote by $\Gamma$ a non-uniform, irreducible, torsion free
lattice in $G$ and hence,
$M=\Gamma\backslash X=\Gamma\backslash G/K$ is a non-compact locally symmetric space
with finite volume. 
From Margulis' famous arithmeticity result it follows that such a $\Gamma$ is always arithmetic 
if $\rank(X)\geq 2$ (\cite{MR1090825,MR776417}). In the rank one case however, it is known
that non-arithmetic lattices exist (\cite{MR932135,MR1090825}). 
Since we treated the rank one case already in  \cite{Ji:yq}, we will restrict ourselves here to 
arithmetic lattices. \\
  
We will now recall some basic facts about the geometry and $L^2$ spectral theory of
arithmetic locally symmetric spaces in order to fix notation. More details can be found e.g. in
\cite{MR0232893,MR1839581,MR1025165}. 
\subsubsection{Langlands decomposition of rational parabolic subgroups}
Since $G= \isom^0(X)$ is a non-compact, semi-simple Lie group with trivial center,
we can find a connected, semi-simple algebraic group $\bG\subset \glnc$ defined over $\Q$ 
such that the groups $G$ and $\bG(\R)^0$ are isomorphic as Lie groups 
(cf. \cite[Proposition 1.14.6]{MR1441541}). \\
A closed subgroup $\bP\subset \bG$ defined over $\Q$ is called {\em rational parabolic subgroup}
if  $\bP$ contains a maximal, connected solvable subgroup of $\bG$.\\
For any rational parabolic subgroup $\bP$ its real locus $P=\bP(\R)$ admits a so-called Langlands decomposition
\begin{equation}
 P = N_{\bP}A_{\bP}M_{\bP}.
\end{equation}
Here $N_{\bP}=\bN_{\bP}(\R)$ denotes the real points of the unipotent radical $\bN_{\bP}$ of $\bP$,
$A_{\bP} = \bS_{\bP}(\R)^0$, where   $\bS_{\bP}$ denotes the maximal $\Q$-split torus
in the center of the Levi quotient $\bL_{\bP} = \bP/\bN_{\bP}$ and hence, $A_{\bP}$
is abelian, and $M_{\bP}$ are the real points of a reductive algebraic group $\bM_{\bP}$ 
defined over $\Q$.
More precisely, this means that the map
\begin{equation}
 P\to N_{\bP}\times A_{\bP}\times M_{\bP},\quad
      g\mapsto \left( n(g), a(g), m(g)\right)
\end{equation}
is a real analytic diffeomorphism.  

Note that such a Langlands decomposition depends on a choice of a maximal compact subgroup
$K$ in $G$ (or equivalently on a base point $x_0\in X$) and $\bM_{\bP}$ is not for all such choices
defined over $\Q$. However, it is known that there exists always a maximal compact subgroup 
$K$ such that the algebraic group $\bM_{\bP}$ is defined over $\Q$. For more details
we refer the reader to the discussion in \cite{MR1906482}.

If we denote by $X_{\bP}$ the {\em boundary symmetric space}
$$ 
 X_{\bP} := M_{\bP}/ K\cap M_{\bP}
$$
we obtain the {\em rational horocyclic decomposition} of $X$:
$$
 X\cong N_{\bP}\times A_{\bP}\times X_{\bP},
$$ 
since the subgroup $P$ acts transitively on the symmetric space $X=G/K$.
More precisely, if we denote by $\tau: M_{\bP}\to X_{\bP}$ the canonical projection, 
we have an analytic diffeomorphism
\begin{equation}\label{rational horocyclic decomposition}
 \mu: N_{\bP}\times A_{\bP}\times X_{\bP} \to X,\,\, (n,a,\tau(m)) \mapsto nam\cdot x_0,
\end{equation}
where $x_0\in X$ denotes a certain base point.\\
Note that the boundary symmetric space $X_{\bP}$ is a Riemannian product of a symmetric
space of non-compact type by possibly a Euclidean space. 
%
\subsubsection{$\Q$-roots and reduction theory}\label{Q-Roots and Reduction Theory}
Let us fix some proper rational parabolic subgroup $\bP$ of $\bG$. We denote in the
following by $\mathfrak{g}, \mathfrak{a}_{\bP}$, and $\mathfrak{n}_{\bP}$ the Lie algebras of the 
real Lie groups $G, A_{\bP}$, and $N_{\bP}$. 
Associated with the pair $(\mathfrak{g}, \mathfrak{a}_{\bP})$ there is  a system 
$\Phi(\mathfrak{g}, \mathfrak{a}_{\bP})$ of  so-called {\em $\Q$-roots} and for each
$\alpha\in \Phi(\mathfrak{g}, \mathfrak{a}_{\bP})$ we have a {\em root space}
$$ 
 \mathfrak{g}_{\alpha} = 
   \{ Z\in \mathfrak{g} : \ad(H)(Y) = \alpha(H)(Y) \mbox{~for all~} H\in  \mathfrak{a}_{\bP} \}.
$$
These root spaces yield a decomposition
$$ 
 \mathfrak{g} = 
    \mathfrak{g}_0 \bigoplus_{\alpha\in \Phi(\mathfrak{g}, \mathfrak{a}_{\bP})}\mathfrak{g}_{\alpha},
$$
where $\mathfrak{g}_0$ is the Lie algebra of $Z(\bS_{\bP}(\R))$, the centralizer of $\bS_{\bP}(\R)$. 
Furthermore, the rational parabolic subgroup $\bP$ defines an ordering of 
$\Phi(\mathfrak{g}, \mathfrak{a}_{\bP})$ such that
$$ 
  \mathfrak{n}_{\bP} = 
     \bigoplus_{\alpha\in \Phi^+(\mathfrak{g}, \mathfrak{a}_{\bP})} \mathfrak{g}_{\alpha},
$$
and the root spaces $\mathfrak{g}_{\alpha}, \mathfrak{g}_{\beta}$ of distinct  positive roots 
$\alpha, \beta\in \Phi^+(\mathfrak{g}, \mathfrak{a}_{\bP})$ are orthogonal with respect to the Killing form:
$$ 
 B(\mathfrak{g}_{\alpha}, \mathfrak{g}_{\beta}) = \{0\}.
$$
We also define
$$ 
 \rho_{\bP} = 
 	\frac{1}{2}\sum_{\alpha\in\Phi^{+}(\mathfrak{g}, \mathfrak{a}_{\bP})}(\dim\mathfrak{g}_{\alpha})\alpha,
$$   
and denote by $\Phi^{++}(\mathfrak{g}, \mathfrak{a}_{\bP})$ the set of simple positive
roots, i.e. the set of all $\alpha\in \Phi^{+}(\mathfrak{g}, \mathfrak{a}_{\bP})$ such that
$\frac{1}{2}\alpha$ is not a root.\\

Let us now discuss reduction theories which describe the structure of fundamental sets for 
$\Gamma$ in terms of Siegel sets associated with rational parabolic subgroups.\\
If we define for $t\in \mathfrak{a}_{\bP}$ 
$$
 A_{\bP,t} = \{ e^{H} \in A_{\bP} : \alpha(H) > \alpha(t) \mbox{~for all~} 
 			\alpha\in \Phi^{++}(\mathfrak{g}, \mathfrak{a}_{\bP}) \},
$$ 
a {\em Siegel set} (associated with $\bP$) is a subset of $X = N_{\bP}\times A_{\bP}\times X_{\bP}$ 
of the form $U\times A_{\bP,t}\times V$ with bounded $U\subset N_{\bP}$
and $V\subset X_{\bP}$.

As it turns out, Siegel sets associated with minimal rational parabolic subgroups are the building
blocks of a fundamental set for $\Gamma$ (for a proof see e.g. \cite[Theorem 13.1]{MR0244260}):
\begin{proposition}\label{non-precise}
 Let $\bP_1,\ldots, \bP_k$ denote representatives of (the finitely many) $\Gamma$ conjugacy classes of
 minimal rational parabolic subgroups. Then there are associated Siegel
 sets ${\cal S}_1,\ldots, {\cal S}_k$ such that
 $F = \bigcup_{j=1}^k {\cal S}_j$
covers a fundamental domain for $\Gamma$ and for any $g\in \bG(\Q)$
the set
$\{\gamma : gF \cap \gamma F\neq\emptyset\}$
is finite, i.e. $F$ is a fundamental set.
\end{proposition}
If we take into account all rational parabolic subgroups of $\bG$, a refined 
(precise) reduction theory yields even a  disjoint decomposition for 
$\Gamma$. 
As Proposition \ref{non-precise} suffices for our purposes, we only refer to
e.g. \cite[Proposition III.2.21]{MR2189882} for further details. 
%
\subsubsection{$\boldsymbol{L^2}$ spectral theory and Eisenstein series}\label{L2 spectral theory}
One knows that the $L^2$ spectrum $\sigma(\Delta_{M,2})$ of the Laplace-Beltrami operator 
$\Delta_{M,2}$  on a non-compact arithmetic locally symmetric space $M$ is the union of a point spectrum and an absolutely continuous spectrum.
The point spectrum consists of a (possibly infinite) sequence of eigenvalues
$$
0=\lambda_0 < \lambda_1\leq \lambda_2\leq \dots
$$
with finite multiplicities such that below any finite number there are only finitely many eigenvalues.
The absolutely continuous spectrum equals $[b,\infty)$ for $b = \min_{\bP} ||\rho_{\bP}||^2$ where
the minimum is taken over all proper rational parabolic subgroups $\bP$.
In what follows, we denote by $L^2_{dis}(M)$ the subspace spanned by all 
eigenfunctions of $\Delta_{M,2}$
and by $L^2_{con}(M)$ the orthogonal complement of  $L^2_{dis}(M)$ in $L^2(M)$.

Generalized eigenfunctions for the absolutely continuous part $\sigma_{ac}(\Delta_{M,2})$ of the $L^2$ spectrum  are given by Eisenstein series.
Therefore, we recall several basic facts about these important functions. Our main reference here is
 \cite{MR0232893}.
\begin{definition}
Let $f$ be a measurable, locally integrable function on $\Gamma\backslash X$. 
The {\em constant term} $f_{\bP}$ of $f$ along some rational parabolic subgroup $\bP$ of 
$\bG$ is defined as
$$
  f_{\bP}(x) = \int_{(\Gamma_{\bP}\cap N_{\bP})\backslash N_{\bP}} f(nx) dn, 
$$
where $\Gamma_{\bP} = \Gamma\cap P$ and the measure $dn$ is normalized such that the  volume of 
$(\Gamma_{\bP}\cap N_{\bP})\backslash N_{\bP}$ equals one.
Note that $(\Gamma_{\bP}\cap N_{\bP})\backslash N_{\bP}$ is always compact if $\bP$ is a rational
parabolic subgroup. \\
 A function  $f$ on $\Gamma\backslash X$  with the property 
 $f_{\bP}=0$ a.e. for all rational parabolic subgroups $\bP\neq \bG$ is called
 {\em cuspidal} and the subspace of cuspidal functions in $L^2(\Gamma\backslash X)$ is 
 denoted by $L^2_{cus}(\Gamma\backslash X)$.
\end{definition}
It is known that 
$$
L^2_{cus}(M) \subset L^2_{dis}(M)
$$ 
and this inclusion is in general strict as the non-zero constant functions are not contained in
$L^2_{cus}(M)$ if $M$ is non-compact.

Let $\bP$ be a rational parabolic subgroup of $\bG$ and $\Gamma_{M_{\bP}}$
the image of $\Gamma_{\bP}=\Gamma\cap P$ under the projection
$N_{\bP}A_{\bP}M_{\bP} \to M_{\bP}$. Then $\Gamma_{M_{\bP}}$ acts discretely on the
boundary symmetric space $X_{\bP}$ and the respective quotient 
$\Gamma_{M_{\bP}}\backslash X_{\bP}$, called {\em boundary locally symmetric space}, has
finite volume. Furthermore, we denote by $\mathfrak{a}_{\bP}^*$ the dual of $\mathfrak{a}_{\bP}$
and put
$$ 
	\mathfrak{a}_{\bP}^{*+} = \{ \lambda \in \mathfrak{a}_{\bP}^* :  
		\langle \lambda, \alpha\rangle > 0  \mbox{~for all~} 
		\alpha \in \Phi^{++}(\mathfrak{g}, \mathfrak{a}_{\bP}) \}.
$$		
For any $\varphi\in L^2_{cus}(\Gamma_{M_{\bP}}\backslash X_{\bP})$ and 
$\Lambda \in \mathfrak{a}_{\bP}^*\otimes\C$ with 
$\re(\Lambda) \in \rho_{\bP} + \mathfrak{a}_{\bP}^{*+}$
we define the {\em (cuspidal) Eisenstein series} $E(\bP|\varphi,\Lambda)$ as follows:
\begin{equation}\label{eisenstein series}
 E(\bP|\varphi,\Lambda:x) = \sum_{\gamma\in \Gamma_{\bP}\backslash \Gamma}
 	e^{(\rho_{\bP}+\Lambda)(H_{\bP}(\gamma x))}\varphi(z_{\bP}(\gamma x)),
\end{equation}
where $\mu(n_{\bP}(x),e^{H_{\bP}(x)}, z_{\bP}(x)) = x$ (cf. (\ref{rational horocyclic decomposition})).
This series converges uniformly for $x$ in compact subsets of $X$ and is holomorphic in 
$\Lambda$ (cf. \cite[Lemma 4.1]{MR0579181}). 
Furthermore,
$E(\bP|\varphi,\Lambda)$ can meromorphically be continued (as a function of $\Lambda$) to
$\mathfrak{a}_{\bP}^*\otimes\C$  (cf. \cite[Chapter 7]{MR0579181} or \cite[Theorem 9]{MR0232893}). 

By definition, the Eisenstein series are $\Gamma$ invariant and hence, they define functions on
$M=\Gamma\backslash X$.

Eisenstein series are in general not contained in $L^2(\Gamma\backslash X)$
but  it is known that they satisfy an eigenequation of the Laplacian 
$\Delta$ on $\Gamma\backslash X$:
\begin{lemma}\label{generalized eigenfunctions}
Let $\varphi\in L^2_{cus}(\Gamma_{M_{\bP}}\backslash X_{\bP})$ be an eigenfunction 
of $\Delta_{\Gamma_{M_{\bP}}\backslash X_{\bP},2}$ with 
respect to some eigenvalue $\nu$. Then we have
for any $\Lambda\in\mathfrak{a}_{\bP}^*\otimes\C$ that is not a pole of 
$E(\bP|\varphi,\Lambda)$ the following:
$$ 
\Delta E(\bP|\varphi,\Lambda) = (\nu + ||\rho_{\bP}||^2
		- \langle \Lambda, \Lambda \rangle ) E(\bP|\varphi,\Lambda),
$$
where $\langle \cdot,\cdot\rangle$ denotes a complex bilinear form on $\mathfrak{a}_{\bP}^*\otimes\C$.
\end{lemma}
For a proof we refer to \cite{MR1906482} or \cite[Lemma 2.5]{Ji:yq}.
%
\section{An upper bound of Eisenstein series and the $\boldsymbol{L^p}$ spectrum of the Laplacian}
The estimates in this section are a precise form of the general philosophy that an automorphic form is bounded by its constant terms.

In order to state these results, we recall that two rational parabolic subgroups
$\bP_1, \bP_2 \subset \bG$ are called {\em associate} ($\bP_1\sim\bP_2$) if there is some 
$g\in \bG(\Q)$ with 
\[
 \Ad(g)\mathfrak{a}_{\bP_1} = \mathfrak{a}_{\bP_2}.
\]
The set of such isomorphisms $\mathfrak{a}_{\bP_1}\to \mathfrak{a}_{\bP_2}$
is denoted by $W(\mathfrak{a}_{\bP_1}, \mathfrak{a}_{\bP_2})$ and
$W(\mathfrak{a}_{\bP}) = W(\mathfrak{a}_{\bP}, \mathfrak{a}_{\bP})$.
For $\Lambda\in\mathfrak{a}_{\bP_1}^*$ and 
$w\in W(\mathfrak{a}_{\bP_1}, \mathfrak{a}_{\bP_2})$
we define $w\Lambda = \Lambda\circ w^{-1} \in \mathfrak{a}_{\bP_2}^*.$

Furthermore, for two rational parabolic subgroups $\bP_1\subset\bP_2$ in $\bG$,
and hence $\mathfrak{a}_{\bP_2}\subset \mathfrak{a}_{\bP_1}$, 
we extend an element $\Lambda\in\mathfrak{a}_{\bP_2}^*$ to an element
in $\mathfrak{a}_{\bP_1}^*$ by defining it to be zero on the orthogonal complement
(w.r.t. the Killing form) of $\mathfrak{a}_{\bP_2}$ in $\mathfrak{a}_{\bP_1}$.
\begin{proposition}\label{upper bound}
Denote by ${\mathcal S}$ a Siegel set associated with a minimal rational parabolic
subgroupf $\bP_0$ of $\bG$ and by $\bP$ a proper rational parabolic subgroup.
If $\varphi \in L^2_{cus}(\Gamma_{M_{\bP}}\backslash X_{\bP})$, the Eisenstein series
$E(\bP | \varphi,\Lambda)$ satisfies the following upper bound for $x\in \mathcal S$:
\begin{equation}
 |E(\bP | \varphi,\Lambda : x)| \leq 
 	C \sum_{\bP'\sim\bP, \atop \bP'\supset\bP_0}\sum_{w\in W(\mathfrak{a}_{\bP},\mathfrak{a}_{\bP'})}
		e^{(w\re\Lambda + \rho_{\bP'})(\log a_{\bP_0}(x))},
\end{equation}
where the constant $C>0$ depends only on $\bP, \varphi, \Lambda$ and $\mathcal S$.
\end{proposition}
\begin{proof}
The key point is to determine the so-called cuspidal data for the Eisenstein series
$E(\bP | \varphi,\Lambda)$ along rational parabolic subgroups $\bP'$.
Then the result basically follows from \cite[Lemma I.4.1]{MR1361168}.

We recall the definition of cuspidal data from \cite{MR1361168} in the notations we use here:
Given an automorphic form $\varphi$ on $\Gamma\backslash X$,
for any rational parabolic subgroup $\bP$, the  constant term $\varphi_{\bP}$ of $\varphi$ along $\bP$ has a cuspidal component  $\varphi_{\bP}^{cusp}$  characterized by the condition that
for every cuspidal function $\psi$ on $\Gamma_{M_\bP}\backslash X_\bP$,
$$\langle \psi, \varphi_{\bP}\rangle = \langle \psi, \varphi_{\bP}^{cusp}\rangle,$$
\cite[p.39]{MR1361168}.\\
The cuspidal component $\varphi_{\bP}^{cusp}(x)$ can be written as a finite sum of functions
of the form $Q(\log a_{\bP}(x)) e^{(\Lambda + \rho_{\bP}) (\log a_{\bP}(x))} \psi(z_{\bP}(x))$,
where $Q$ is a polynomial and $\psi$ a cuspidal form on $\Gamma_{M_{\bP}}\backslash X_{\bP}$.
The finite triples $(Q, \Lambda, \psi)$ are called the cuspidal data of $\varphi_{\bP}^{cusp}$,
or the cuspidal data of $\varphi$ along $\bP$, cf. \cite[p.44]{MR1361168}.
Note that in \cite{MR1361168} representations $\pi$ are used instead of characters $\Lambda$ 
(or rather functionals in $\mathfrak{a}_{\bP}^*\otimes\C$). This is the reason for the shift
by $+\rho_{\bP}$ appearing here (cf. also \cite[Chapter VII.1]{0993.22001}).\\

Let us now determine the cuspidal data for the Eisenstein series
$E(\bP | \varphi,\Lambda)$ along some rational parabolic subgroup $\bP'$.

If $\bP, \bP'$ are not associate, the constant term $E_{\bP'}(\bP | \varphi,\Lambda)$ is orthogonal to all cusp forms on $\Gamma_{M_{\bP}}\backslash X_{\bP}$, cf. \cite[Lemma 39]{MR0232893}
or \cite[p.86]{MR643242} and hence, in this case, the set of cupidal data is empty.\\
If on the other hand $\bP$ and $\bP'$ are associate, 
we have the following formula for the constant term of 
$E(\bP | \varphi,\Lambda)$ along $\bP'$:
\begin{equation}
 E_{\bP'}(\bP | \varphi,\Lambda : x) =
 	\sum_{w\in W(\mathfrak{a}_{\bP},\mathfrak{a}_{\bP'})}
	     e^{(w\Lambda + \rho_{\bP'})(\log(a_{\bP'}(x))}\cdot
	     (c_{cus}(\bP' : \bP : w : \Lambda)\varphi)(z_{\bP'}(x)),
\end{equation}
where $c_{cus}(\bP' : \bP : w : \Lambda)$ denotes the intertwining operator from the space
of cusp forms on $\Gamma_{M_{\bP}}\backslash X_{\bP}$ to the space of cusp forms on
$\Gamma_{M_{\bP'}}\backslash X_{\bP'}$, 
cf. \cite[Theorem 5]{MR0232893} or \cite[p.86]{MR643242}. Note that in 
the last mentioned book a slightly different definition of the ``constant term'' of an automorphic form is used,  see \cite[p.79]{MR643242} for the definition used therein.\\
Hence, for associated $\bP, \bP'$ the set of cuspidal data of $E(\bP | \varphi,\Lambda)$ along $\bP'$
consists of the triples
\begin{equation}
 (Q=1, w\Lambda, c_{cus}(\bP' : \bP : w : \Lambda)\varphi), \qquad 
 						w\in W(\mathfrak{a}_{\bP},\mathfrak{a}_{\bP'}).
\end{equation}
From  \cite[Lemma I.4.1]{MR1361168} it now follows that there is some $C>0$ such that 
for $x\in S$ we have
\begin{eqnarray*}
 |E(\bP | \varphi,\Lambda : x)| &\leq& 
 	C\sum_{\bP'\sim\bP, \atop \bP'\supset \bP_0}\sum_{w\in W(\mathfrak{a}_{\bP},\mathfrak{a}_{\bP'})}
		e^{(w\re\Lambda + \rho_{\bP'})(\log a_{\bP_0}(x))}\, .
\end{eqnarray*}
Note that in \cite[Lemma I.4.1]{MR1361168} the summation considers only standard parabolic subgroups, which are defined as rational parabolic subroups $\bP$ with $\bP\supset \bP_0$.
In our situation here, the subgroups $\bP'$ associated with $\bP$ are the only standard parabolic subgroups which lead to non-empty cuspidal data. 
\end{proof}
Recall that all minimal rational parabolic subgroups are conjugate under $\bG(\Q)$. Furthermore,
if $\bP$ is a minimal rational parabolic subgroup, the boundary locally symmetric space 
$\Gamma_{M_{\bP}}\backslash X_{\bP}$  is compact and hence any $L^2$ eigenfunction
$\varphi$ on $\Gamma_{M_{\bP}}\backslash X_{\bP}$ is cuspidal, i.e 
$L^2_{cus}(\Gamma_{M_{\bP}}\backslash X_{\bP}) = L^2(\Gamma_{M_{\bP}}\backslash X_{\bP})$.
Therefore Proposition \ref{upper bound} simplifies in the case where $\bP$ is a minimal rational parabolic subgroup:
\begin{corollary}\label{corollary eisenstein upper bound minimal} 
Let $\bP_0, \bP$ denote minimal rational parabolic subgroups of $\bG$,
${\mathcal S}$ a Siegel set associated with $\bP_0$, and  $E(\bP | \varphi,\Lambda)$ an Eisenstein series associated with $\bP$.
Then there exists a constant $C>0$ such that for all $x\in{\mathcal S}$
\begin{eqnarray*}
 |E(\bP | \varphi,\Lambda : x)|  &\leq& 
 	C\sum_{w\in W(\mathfrak{a}_{\bP},\mathfrak{a}_{\bP_0})}
		e^{(w\re\Lambda + \rho_{\bP_0})(\log a_{\bP_0}(x))}.
\end{eqnarray*}
\end{corollary}
%

\subsection{$\boldsymbol{L^p}$ spectrum} \label{Lp spectrum}
%
Let us define for a minimal parabolic subgroup $\bP_0$ and a rational parabolic subgroup $\bP$ of $\bG$ the set
\begin{multline}
 C({\bP}) = C(\bP,\bP_0) = \\	\bigcap_{\substack{w\in W(\mathfrak{a}_{\bP},\mathfrak{a}_{\bP'}), \\
 																		   \bP' \sim \bP,\, \bP'\supset \bP_0} } 
 \left\{ \Lambda : w\re\Lambda (\log a_{\bP'}) < \left(\frac{2}{p}-1\right)  \rho_{\bP'}(\log a_{\bP'}), a_{\bP'}\in A_{\bP',0} \right\}.
\end{multline}
Note that this is a finite intersection since $W(\mathfrak{a}_{\bP},\mathfrak{a}_{\bP'})$ is finite and since there are only finitely many $\bP' \sim \bP, \bP'\supset \bP_0$. Furthermore, the set is non-empty since it contains $\Lambda = 0$.\\
If $\bP = \bP_0$ has $\Q$-rank one, i.e. $\dim A_{\bP} = 1$, the set simplifies to
$$
 C({\bP}) = 
    \left\{ \Lambda : |\re\Lambda (\log a_{\bP}) | < \left(\frac{2}{p}-1\right) \rho_{\bP}(\log a_{\bP}), a_{\bP}\in A_{\bP,0} \right\},
$$
since $W(\mathfrak{a}_{\bP},\mathfrak{a}_{\bP}) = \{ -\id, \id \}$.
Therefore, the results below are natural extensions of the results in \cite{Ji:yq}, where rank one locally symmetric spaces are treated, to the higher rank case.\\ 
Recall that all proper rational parabolic subgroups are minimal if the locally symmetric space
$M=\Gamma\backslash X$ has $\Q$-rank one. Furthermore, two minimal parabolic subgroups $\bP, \bP'$ are always 
conjugate under $\bG(\Q)$ and hence the group  $W(\mathfrak{a}_{\bP'},\mathfrak{a}_{\bP})$ can be identified
with the Weyl group $W(\mathfrak{a}_{\bP})$. Using this identification, we have
\begin{lemma}\label{convex subset}
If  $\bP = \bP_0$ denotes a minimal parabolic subgroup of $\bG$ it follows
$$
 \mathrm{int}\, \mathrm{conv} \left\{\left(\frac{2}{p} - 1\right) w\rho_{\bP}: w \in W(\mathfrak{a}_{\bP})  \right\} \subset C(\bP),
$$ 
where $\mathrm{conv}$ denotes convex hull.
\end{lemma}
\begin{proof}
If 
$$
\Lambda \in   \mathrm{int}\, \mathrm{conv} \left\{\left(\frac{2}{p} - 1\right) w\rho_{\bP}: w \in W(\mathfrak{a}_{\bP})  \right\},
$$
and $W(\mathfrak{a}_{\bP} ) = \{w_{j} : j \}$
then
$$
\Lambda = \left(\frac{2}{p} - 1\right) \sum_{j} r_{j}w_{j}\rho_{\bP}
$$
for some  $0 < r_{j} < 1$ with $\sum_{j} r_{j} = 1$.\\
Let  now $w_1 \in W(\mathfrak{a}_{\bP})$. Then
$$
w_1\Lambda = \left(\frac{2}{p} - 1\right) \sum_{j} r_{j}w_1w_{j}\rho_{\bP}.
$$
Furthermore, we have for $a_{\bP} \in A_{\bP,0}$ and any $j$ with 
$w_1w_{j}\rho_{\bP} \neq \rho_{\bP}$
$$
 w_1w_{j}\rho_{\bP}(\log(a_{\bP} )) < \rho_{\bP}(\log(a_{\bP} ))
$$
since for any positive root $\alpha \in \mathfrak{a}_{\bP}^*$ the functional
$w_1w_{j}\alpha$ is either negative on  $\mathfrak{a}_{\bP,0}$ or a positive root $\in \mathfrak{a}_{\bP}^*$.\\
For $0 < r_{j} < 1$ with $\sum_{j} r_{j} = 1$ it then follows
\begin{equation*}
w_1\Lambda(\log a_{{\bP}}) = \left(\frac{2}{p} - 1\right) \sum_{j} r_{j} w_1w_{j}\rho_{\bP}(\log a_{{\bP}}) 
  < \left(\frac{2}{p} - 1\right)  \rho_{\bP}(\log a_{{\bP}})
\end{equation*}
and therefore $\Lambda \in C(\bP)$.
\end{proof}
\begin{lemma}\label{lemma estimate siegel set}
Let $\bP$ denote a rational parabolic subgroup, $E(\bP | \varphi, \Lambda)$
an associated Eisenstein series on $M=\Gamma\backslash X$ with 
$\varphi \in L^2_{cus}(\Gamma_{M_P}\backslash X_P)$, ${\mathcal S}$
a Siegel set associated with a minimal rational parabolic subgroup $\bP_0$, and $p\in [1,2)$. 
If 
$$ 
	\Lambda \in  C({\bP})
$$
and, if $\Lambda$ is not a pole of $E(\bP | \varphi, \Lambda)$, we have
$E(\bP | \varphi, \Lambda) \in L^p({\mathcal S})$.
\end{lemma}
\begin{proof}
Recall that the Riemannian measure on $X$ with respect to the horocyclic
decomposition $X = N_{\bP_0}\times A_{\bP_0}\times X_{\bP_0}$ is given by
$dvol_X = h(n_{\bP_0},z_{\bP_0})e^{-2\rho_{\bP_0}(\log a_{\bP_0})}dndzda$, where $h$ is smooth on
$N_{\bP_0}\times X_{\bP_0}$, cf. \cite[Proposition 1.6]{MR0338456} or
\cite[Proposition 4.3]{MR0387496}. Therefore, we obtain
for $p\in [1,2)$ with Proposition \ref{upper bound}
\begin{eqnarray*}
 \int_{{\mathcal S}} |E(\bP | \varphi,\Lambda)(x)|^p dvol_X &\leq& 
 	  C\sum_{\bP'\sim\bP,\atop \bP'\supset\bP_0}
	  \sum_{w\in W(\mathfrak{a}_{\bP},\mathfrak{a}_{\bP'})}\int_{A_{\bP_0, t}} 
		 e^{(pw\re\Lambda + p\rho_{\bP'} - 2\rho_{\bP_0})(\log a_{\bP_0}(x))} da.
\end{eqnarray*}
We note that for any value $t$, $A_{\bP_0, t}$ is a shift of the positive cone
$A_{\bP_0, 0}$,  and we can assume $t=0$ for simplicity,
and hence, $E(\bP | \varphi, \Lambda) \in L^p(S)$ if  
\begin{equation}\label{inequality1}
(pw\re\Lambda + p\rho_{\bP'} - 2\rho_{\bP_0})(\log a_{\bP_0}) < 0
\end{equation}
for all $a_{\bP_0}\in A_{\bP_0, t}, w\in W(\mathfrak{a}_{\bP},\mathfrak{a}_{\bP'})$, and all
rational parabolic subgroups $\bP'\sim\bP$ with $\bP'\supset\bP_0$.

As by definition $w\re\Lambda  = \rho_{\bP'} = 0$ on $\mathfrak{a}_{\bP'}^{\perp}$ (the orthogonal complement of $\mathfrak{a}_{\bP'}$ in $\mathfrak{a}_{\bP_0}$) and since
the restriction of $\rho_{\bP_0}$ to $\mathfrak{a}_{\bP'}$ coincides with $\rho_{\bP'}$ 
(cf. \cite[Lemma 1.14.1]{MR0231321}),
inequality (\ref{inequality1}) is equivalent to
\begin{eqnarray}\label{second inequality}
 (pw\re\Lambda + (p-2)\rho_{\bP'})(\log a_{\bP'}) & < & 0, 
\end{eqnarray}
for all $a_{\bP'}\in A_{\bP',t} = A_{\bP'}\cap A_{\bP_0,t}$ and $w\in W(\mathfrak{a}_{\bP},\mathfrak{a}_{\bP'})$. 
This proves the claim. 
\end{proof}
We are now prepared to prove the following result.
\begin{theorem}\label{thm general parabolic}
Let $\bP$ denote a rational parabolic subgroup of $\bG$, $E(\bP | \varphi, \Lambda)$
an associated Eisenstein series on $M=\Gamma\backslash X$ 
with $\varphi \in L^2_{cus}(\Gamma_{M_P}\backslash X_P)$
and $p\in [1,2)$.
If
\begin{eqnarray}\label{condition lambda X}
	\Lambda \in  C({\bP})
\end{eqnarray}
and if $\Lambda$ is not a pole, we have 
$E(\bP | \varphi, \Lambda) \in L^p(M)$.
\end{theorem}
\begin{proof}
Let $\bP_j, j=0,\ldots,k,$ denote representatives of $\Gamma$ conjugacy classes of 
minimal rational parabolic subgroups.
From Proposition \ref{non-precise} it follows that there are Siegel sets
${\cal S}_j$ associated with $\bP_j, j=0,\ldots,k,$ such that the union
of these Siegel sets covers a fundamental domain for $\Gamma$.
From Lemma \ref{lemma estimate siegel set} it follows
$E(\bP | \varphi, \Lambda) \in L^p({\cal S}_j)$ if
\begin{eqnarray*}
	\Lambda \in  C({\bP})
\end{eqnarray*}
\end{proof}
\begin{corollary}\label{eigenvalues minimal}
Let $\bP$ denote a rational parabolic subgroup of $\bG$, $E(\bP | \varphi, \Lambda)$
an associated Eisenstein series on $M=\Gamma\backslash X$ 
with $\varphi \in L^2_{cus}(\Gamma_{M_{\bP}}\backslash X_{\bP})$ and $p\in (1,2)$.
If
\begin{eqnarray*}
	\Lambda \in  C({\bP})
\end{eqnarray*}
and if $\Lambda$ is not a pole, $E(\bP | \varphi, \Lambda)$ is an eigenfunction
of $\DMp$ with eigenvalue $\nu + ||\rho_{\bP}||^2 - \langle \Lambda,\Lambda\rangle$, where
$\nu$ is the eigenvalue of the $L^2$ Laplacian on the boundary locally symmetric space 
$\Gamma_{M_{\bP}}\backslash X_{\bP}$ for $\varphi$.
\end{corollary}
\begin{proof}
This follows from Lemma \ref{lemma Lp eigenfunctions} and 
Lemma \ref{generalized eigenfunctions} together with the preceeding theorem.
\end{proof}
These results resemble the behavior of spherical functions 
$$
\varphi_{\lambda}(gK) = \int_K \e^{(i\lambda + \rho)(A(kg))} dk,\qquad \lambda\in\mathfrak{a}^*_{\C}
$$
on (globally) symmetric spaces $X=G/K$ of non-compact type: If $C_X(\rho)$ denotes 
the convex hull of the points $s\rho\in\mathfrak{a}^*, s\in W$, where $W$ 
denotes the Weyl group, for any $p>2$ and any 
$\lambda \in \mathfrak{a}^* + i(1-\frac{2}{p})C_X(\rho)$ the spherical function $\varphi_{\lambda}$ is contained in $L^p(X)$, cf. \cite[Proposition 2.2]{MR1016445} (note that the roles of $p>2$ and $p<2$
are interchanged). \\

If we define for $p\in [1,\infty)$  the parabolic region
\begin{equation}
 {\cal P}_{M,p}(\rho_{\bP}) = \left\{ ||\rho_{\bP}||^2 - z^2 : 
 			z\in\C, |\re z|\leq ||\rho_{\bP}||\cdot |\frac{2}{p}-1|\right\}\subset \C,
\end{equation} 
we obtain as in \cite[Proposition 2.2]{MR1016445} or \cite[Corollary 3.6]{Ji:yq}
the following:
\begin{theorem}\label{theorem eigenvalues}
Let $\bP$ denote a minimal rational parabolic subgroup of $\bG$.
\begin{itemize}
\item[\textup{(a)}]  If $p\in (1,2)$, there exists a discrete set 
		$B\subset {\cal P}_{M,p}(\rho_{\bP})$ such that all points in the interior of 
		${\cal P}_{M,p}(\rho_{\bP}) \setminus B$ are eigenvalues of $\DMp$.
\item[\textup{(b)}] If $p\in [2,\infty)$, the point spectrum $\sigma_{pt}(\DMp)$ is a 
		discrete subset of $[0,\infty)$.		
\end{itemize}
Furthermore, ${\cal P}_{M,p}(\rho_{\bP}) \subset \sigma(\DMp)$ for all $p\in[1,\infty)$.
\end{theorem}
Recall that two minimal rational parabolic subgroups $\bP,\bP'$ are conjugate under 
$\bG(\Q)$ and in particular, $||\rho_{\bP}|| = ||\rho_{\bP'}||$. Therefore we have
${\cal P}_{M,p}(\rho_{\bP}) = {\cal P}_{M,p}(\rho_{\bP'})$.
\begin{proof}
Let first $p\in (1,2)$ and let $E(\bP | \varphi, \Lambda)$ denote
an Eisenstein series associated with $\bP$ on $M=\Gamma\backslash X$ 
with $\varphi= const \in L^2_{cus}(\Gamma_{M_{\bP}}\backslash X_{\bP})$. Note that such a choice is always possible as -- due to  the minimality of $\bP$ -- 
the boundary locally symmetric space $\Gamma_{M_{\bP}}\backslash X_{\bP}$ is compact
and hence, $L^2_{cus}(\Gamma_{M_{\bP}}\backslash X_{\bP}) =
 L^2(\Gamma_{M_{\bP}}\backslash X_{\bP})$.\\
From Corollary \ref{eigenvalues minimal} and Lemma \ref{convex subset} it follows that for $\Lambda$ with
\begin{equation}
0 < \re\Lambda (\log a_{\bP}) < \left(\frac{2}{p}-1\right) \rho_{\bP}(\log a_{\bP}),\quad a_{\bP}\in A_{\bP,0},
\end{equation}
or equivalently if
$$
\Lambda \in \left\{ z\Lambda_0: z\in\C, 0 < \re(z) < \left(\frac{2}{p}-1\right) || \rho_{\bP} ||\right\}
$$
for some $\Lambda_0\in \mathfrak{a}_{\bP}^*$ with $|| \Lambda_0 || = 1$,
all the points 
$||\rho_{\bP}||^2 - \langle \Lambda, \Lambda\rangle$
are contained in the point spectrum of $\DMp$ if $\Lambda$ is not a pole. 
Therefore there exists a discrete set $B$ such that all points in the interior of
${\cal P}_{M,p}(\rho_{\bP}) \setminus B$
are eigenvalues of $\DMp$ (recall that the poles of Eisenstein series lie along affine hyperplanes). Since the spectrum is a closed subset of $\C$, we obtain
${\cal P}_{M,p}(\rho_{\bP}) \subset \sigma(\DMp)$ for $p\in (1,2)$ and by duality
for $p\in(2,\infty)$. In the case $p=2$, we have 
${\cal P}_{M,2}(\rho_{\bP}) = [||\rho_{\bP}||^2,\infty)$, and it is well known that this set
belongs to the $L^2$ spectrum of $M$. This completes the proof of (a) and the last statement.\\
To prove (b), assume $p>2$ and $\DMp \psi = \lambda \psi$ for some $\psi\in \dom(\DMp)$.
Since $vol(M)<\infty$ it follows from Lemma \ref{domains} that $\psi$ is an eigenfunction of
$\DMq$ with eigenvalue $\lambda$ if $1<q\leq p$. Since it is known that the point spectrum of 
$\Delta_{M,2}$
is a discrete subset of $[0,\infty)$, the proof is complete.
\end{proof}
An application of Theorem \ref{thm general parabolic} yields similarly for non-minimal rational
parabolic subgroups the existence of an open subset of $\C$ that consists only of eigenvalues
of $\DMp$. However, an explicit description of this set seems to be more complicated
than in the case of minimal rational parabolic subgroups.\\

\begin{figure}[htb]
  \centering
  \includegraphics[scale=0.4]{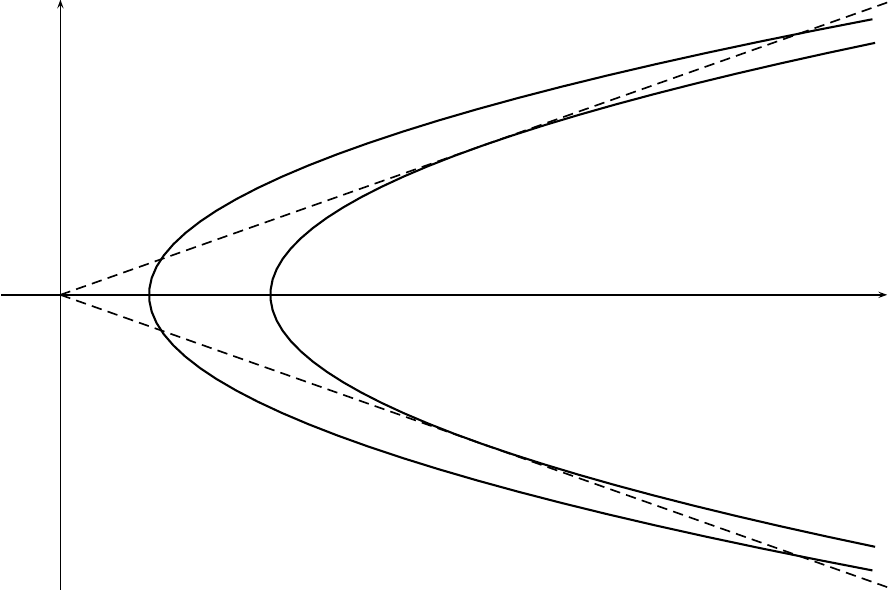}
   \caption{The parabolic region  ${\cal P}_{M,p}(\rho_{\bP})$ tangent to the sector
   		$\Sigma_p$ defined in Section \ref{heat semigroup} and the parabolic 
		region $ {\cal P}_{M,p}'$.} 
   \label{parabolic regions}
\end{figure}

On the other hand, an ``upper bound'' for the $L^p$ spectrum
of $M$ follows from a result by Taylor:
\begin{proposition}
Let $M=\Gamma\backslash X$ denote a non-compact arithmetic locally symmetric space 
and denote by $\bP$ a minimal rational parabolic subgroup.
For $p\in (1,\infty)$
we have
$$ 
 \{\lambda_0,\ldots, \lambda_r\}\cup {\cal P}_{M,p}(\rho_{\bP})\subset \sigma(\DMp) \subset 
 		\{\lambda_0,\ldots, \lambda_r\}\cup {\cal P}_{M,p}',
$$
where 
$$
 {\cal P}_{M,p}' = \left\{ b - z^2 :  z\in\C, |\re z|\leq ||\rho||\cdot |\frac{2}{p}-1|\right\}
$$
with $b = \inf \sigma_{ac}(\Delta_{M,2}) = \inf_{\bP\neq \bG} ||\rho_{\bP}||^2 > 0,$
and where $\lambda_0 = 0, \lambda_1,\ldots, \lambda_r$ are eigenvalues of $\Delta_{M,2}$
with finite multiplicity.
\end{proposition}
\begin{proof}
The second inclusion follows from Taylor's Proposition 3.3 in \cite{MR1016445} (see also
Proposition \ref{taylor upper} in this paper). 
Taking into account Theorem \ref{theorem eigenvalues}, it remains to show
$\lambda_j \in \sigma(\DMp), j=0,\ldots, r$. For $p\leq 2$ it follows similarly as
in the proof of Theorem \ref{theorem eigenvalues} that each $\lambda_j$ is
an eigenvalue of $\DMp$.
By duality, it then follows $\lambda_j\in \sigma(\DMp)$ for any $p>2$ and the
proof is complete.
\end{proof}
If $X=G/K$ denotes a symmetric space of non-compact type and $M=\Gamma\backslash X$
for a non-uniform arithmetic $\Gamma$,
in the case where $\qrank(\Gamma) = \rank(X)$, we have for any minimal rational parabolic
subgroup $\bP\subset\bG$
$$
\dim(A_{\bP}) = \dim(A),
$$
where $A$ denotes the abelian subgroup in the Iwasawa decomposition $G=KAN$.
Hence, in this situation, we can conclude $||\rho|| = ||\rho_{\bP}||$ and therefore
${\cal P}_{M,p}'$ is a shift of ${\cal P}_{M,p}(\rho_{\bP})$, more 
precisely ${\cal P}_{M,p}' = (b-||\rho_{\bP}||^2) + {\cal P}_{M,p}(\rho_{\bP})$.
See Figure \ref{parabolic regions}.

%

%

\begin{figure}[htb]
  \centering
  \includegraphics[scale=0.7]{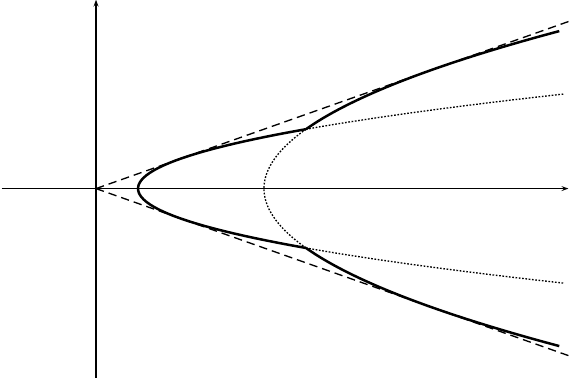}
   \caption{Conjecture about the $L^p$ spectrum with sector
   				$\Sigma_p$.} 
   \label{conjecture figure}
\end{figure}

We conclude this section with a conjecture about the precise form of the $L^p$ spectrum.
\begin{conjecture}[cf. Figure \ref{conjecture figure}] \label{conjecture}
Let $M=\Gamma\backslash X$ denote a non-compact locally symmetric space 
with arithmetic fundamental group $\Gamma$ and denote by $0=\lambda_0,\ldots,\lambda_r$
the eigenvalues of $\Delta_{M,2}$ below the absolutely continuous spectrum
$\sigma_{ac}(\Delta_{M,2})$. Then, for any $p\in(1,\infty)$, 
each  proper rational 
parabolic subgroup $\bP$ defines a parabolic 
region ${\cal P}_{\bP,p}$ tangent to the boundary of the sector 
$\Sigma_p = \left\{ z\in \C\setminus\{0\} : |\arg(z)| \leq \arctan\frac{|p-2|}{2\sqrt{p-1}} \right\}$
with apex at the point $z(\bP)=\frac{4||\rho_{\bP}||^2}{p}\left(1- \frac{1}{p}\right)$
such that the $L^p$ spectrum $\sigma(\DMp)$ of $M$ coincides precisely
with the union of these parabolic regions ${\cal P}_{\bP,p}$
and the $L^2$ eigenvalues $\lambda_0,\ldots,\lambda_r$:
$$
 \sigma(\DMp) = \{\lambda_0,\ldots,\lambda_r\}\cup \bigcup_{\bP} {\cal P}_{\bP,p}.
$$
Furthermore, if $p\in (1,2)$, there is a discrete subset $B\subset \C$ such that 
all points in 
$\left(\{\lambda_0,\ldots,\lambda_r\}\cup \bigcup {\cal P}_{\bP,p}\right)\setminus B$
are eigenvalues of $\DMp$.
\end{conjecture}
We actually assume that, similarly as in the case of minimal rational parbolic subgroups, each
{\em non-cuspidal} Eisenstein series $E(\bP|\varphi_0,\Lambda)$ with $\varphi_0=const.$ defines a parabolic region ${\cal P}_{\bP,p}$. For fixed $\bP$, it seems to be plausible that for any
$\varphi\in L^2(\Gamma_{M_{\bP}}\backslash X_{\bP})$ the parabolic region  defined by 
$E(\bP|\varphi,\Lambda)$ should be contained in ${\cal P}_{\bP,p}$. 
More precisely, for any 
$\varphi\in L^2(\Gamma_{M_{\bP}}\backslash X_{\bP})$ 
the parabolic region defined by $E(\bP|\varphi,\Lambda)$ should be a shift (by the eigenvalue corresponding to $\varphi$) of  the parabolic region defined by $E(\bP|\varphi_0,\Lambda)$.\\
In the case $p=2$ each parabolic region ${\cal P}_{\bP,p}$ 
of course needs to degenerate to a subset of $\R$ as $\Delta_{M,2}$ is self-adjoint. 
It is very plausible that this subset is of the form $[||\rho_{\bP}||^2,\infty)$.
Going from $L^2$ spectrum to $L^p$ spectrum, $p\neq 2$ would then be like a butterfly opening its wings.\\
If Conjecture \ref{conjecture} is true, the $L^p$ spectrum (as set), for any $p\neq 2$, 
contains more information about $M$ than the $L^2$ spectrum. For example, if the $L^p$
spectrum containes $n$ different parabolic regions (i.e. $n$ parabolic regions tangent to $\Sigma_p$ with pairwise different apex), there must be $n$ different rational parabolic
subgroups $\bP$ such that $||\rho_{\bP}||$ are pairwise different numbers.\\

Conjecture \ref{conjecture} is known to be true for non-compact locally symmetric spaces
$M=\Gamma\backslash X$ with finite volume if $X$ is a rank one symmetric space
of non-compact type \cite{MR937635,Ji:yq}, for Hilbert modular varieties $M$ \cite{Ji:yq},
and by the result in \cite{Weber:2008ve} for any Riemannian product of these spaces.

\section{Heat dynamics}
In this section we use our results of the previous sections to show that the
$L^p$ heat semigroup has a certain chaotic behavior if $p\in(1,2)$ whereas such 
a chaotic behavior cannot occur if $p\geq 2$. \\
Similar results have been proven in the context of globally symmetric spaces of non-compact type
 \cite{Ji:nr} and for rank one locally symmetric spaces (which need not neccessarily be arithmetic)
 in \cite{Ji:yq}.

\subsection{Chaotic semigroups}

There are many different definitions of chaos. We will use the following one which is
basically an adaption of Devaney's definition \cite{MR1046376} to the setting of strongly
continuous semigroups (cf. \cite{MR1468101}).
\begin{definition}
 A strongly continuous semigroup $T(t)$ on a Banach space ${\cal B}$ is called {\em chaotic}
 if the following two conditions are satisfied:
   \begin{itemize}
    \item[\textup{(i)}] There exists an $f\in {\cal B}$ such that its orbit 
    			      $\{T(t)f : t\geq 0 \}$ is dense in ${\cal B}$, i.e. $T(t)$ is {\em hypercyclic}.	
    \item[\textup{(ii)}] The set of periodic points
    			        $\{ f\in {\cal B} : \exists t>0 \mbox{~such that~} T(t)f=f \}$	 is dense in ${\cal B}$.		   
   \end{itemize}
\end{definition}
\begin{remark}\label{remark1}
 \begin{itemize}
  \item[\textup{(1)}] As with $\{T(t)f : t\geq 0 \}$  also the set $\{T(q)f : q\in\Q_{\geq 0} \}$
  				is dense, ${\cal B}$ is necessarily separable.
  \item[\textup{(2)}] The orbit of any point $T(t)f$ in a dense orbit $\{T(t)f : t\geq 0 \}$ is 
  				again dense in ${\cal B}$. Hence, the set of 
				points with a dense orbit is a dense subset of ${\cal B}$ or empty.
    \item[\textup{(3)}] If $A$ is the generator of a hypercyclic semigroup, its dual operator
    			       $A'$ has empty point spectrum \cite[Theorem 3.3]{MR1468101}.	
  \end{itemize}
\end{remark}
A sufficient condition for a strongly continuous semigroup to be chaotic in terms of spectral properties
of its generator was given by Desch, Schappacher, and Webb: 
\begin{theorem}[\cite{MR1468101}] \label{thm dsw}
Let $T(t)$ denote a strongly continuous semigroup on a separable 
 Banach space ${\cal B}$ with generator $A$ and let  $\Omega$ denote an open, connected
 subset of $\C$ with  $\Omega\subset \sigma_{pt}(A)$  (the point spectrum of $A$). 
 Assume that there is  a function $F: \Omega\to {\cal B}$ such that
 \begin{itemize}
  \item[\textup{(i)}] $\Omega \cap i\R \neq \emptyset$.
  \item[\textup{(ii)}] $F(\lambda) \in \ker(A-\lambda)$ for all $\lambda \in \Omega$.
  \item[\textup{(iii)}] For all $\phi \in {\cal B'}$ in the dual space of ${\cal B}$, the mapping
  				$F_{\phi}:\Omega\to \C,\, \lambda\mapsto \phi\circ F $
				is analytic. 
				Furthermore, if for some $\phi \in {\cal B'}$ we have $F_{\phi}=0$
				then already $\phi = 0$ holds.
 \end{itemize}
 Then the semigroup $T(t)$ is chaotic.
\end{theorem}
In  \cite{MR1468101} it was also required that the elements $F(\lambda)$, $\lambda \in \Omega$, are non-zero but as remarked in \cite{MR2128736} this assumption is redundant. \\

In the theory of dynamical systems chaotic semigroups are highly unwanted because of their
difficult dynamics. Not much more appreciated are so called subspace chaotic semigroups:
\begin{definition} 
  A strongly continuous semigroup $T(t)$ on a Banach space ${\cal B}$ is called {\em subspace 
  chaotic} if there is a closed, $T(t)$ invariant subspace ${\cal V}\neq \{0\}$ of  ${\cal B}$ such that
  the restriction $T(t)|_{\cal V}$ is a chaotic semigroup on ${\cal V}$.
\end{definition}
\noindent Because of Remark \ref{remark1} such a subspace is always infinite dimensional.\\

Banasiak and Moszy\'nski showed that a subset of the conditions in Theorem \ref{thm dsw} yields a sufficient  condition for subspace chaos:
\begin{theorem}\textup{(\cite[Criterion 3.3]{MR2128736}).}\label{thm ban}
Let $T(t)$ denote a strongly continuous semigroup on a separable 
 Banach space ${\cal B}$ with generator $A$. Assume, there is an open, connected subset
 $\Omega\subset \C$ and a function $F: \Omega\to {\cal B}, F\neq 0,$ such that
 \begin{itemize}
  \item[\textup{(i)}] $\Omega \cap i\R \neq \emptyset$.
  \item[\textup{(ii)}] $F(\lambda) \in \ker(A-\lambda)$ for all $\lambda \in \Omega$.
  \item[\textup{(iii)}] For all $\phi \in {\cal B'}$, the mapping
  				$F_{\phi}:\Omega\to \C,\, \lambda\mapsto \phi\circ F $
				is analytic. 
 \end{itemize}
 Then the semigroup $T(t)$ is subspace chaotic.\\
 Furthermore, the restriction of $T(t)$ to the $T(t)$ invariant subspace 
 ${\cal V} = \overline{\spa}F(\Omega)$ is chaotic.	  
\end{theorem}
Note that it is not required $\Omega\subset \sigma_{pt}(A)$ here, i.e. either $F(\lambda)$ is an eigenvector  or $F(\lambda)=0$.  But, as explained in \cite{MR2128736}, the assumption $\Omega\subset \C$ is not really weaker.
\subsection{$\boldsymbol{L^p}$ heat dynamics on locally symmetric spaces}
%
\begin{theorem}\label{thm dynamics}
Let $M=\Gamma\backslash X$ denote a non-compact locally symmetric space with
arithmetic fundamental group $\Gamma$.
\begin{itemize}
 \item[\textup{(a)}]
 If $p\in (1,2)$ there is a constant $c_p>0$ such that for any $c>c_p$ the semigroup
 $$
  e^{-t(\DMp -c)}: L^p(M) \to L^p(M)
 $$
 is subspace chaotic.
 \item[\textup{(b)}]
 For any $p\geq 2$ and $c\in\R$ the semigroup $e^{-t(\DMp -c)}: L^p(M) \to L^p(M)$ 
 is not subspace chaotic.
\end{itemize}
\end{theorem}
\begin{proof}
For the proof of part (a), we will check the conditions of Theorem  \ref{thm ban}.
If $p\in (1,2)$ and if $\bP$ denotes a minimal rational parabolic subgroup
of $\bG$, the interior of 
$$
 ({\cal P}_{M,p}(\rho_{\bP})\setminus B)\cap\{ z\in\C : \im(z) <0\}
$$ 
consists completely of eigenvalues (cf. Theorem \ref{theorem eigenvalues}), and the apex of 
${\cal P}_{M,p}(\rho_{\bP})$
is at the point  
$$
 c_p = ||\rho_{\bP}||^2 - ||\rho_{\bP}||^2 \cdot \left(\frac{2}{p} - 1\right)^2 
         = \frac{4||\rho_{\bP}||^2}{p}\left(1 - \frac{1}{p}\right).
$$
Hence, the point spectrum of $(\DMp - c)$ intersects the imaginary axis for any $c>c_p$.
We assume in the following $c > c_p$ and denote by $\Omega$ the interior of the set
$$
  \left({\cal P}_{M,p}(\rho_{\bP})\setminus B - c\right)\cap\{ z\in\C : \im(z) <0\}.
$$
Then, if the usual analytic branch of the square root is chosen, 
$\Omega$ is mapped (analytically) by 
$h(z) = i ||\rho_{\bP}||^{-1}\sqrt{z+c-||\rho_{\bP}||^2}$
onto the strip
$$
 \left\{z\in\C : \im(z)>0, 0<\re(z)< \left(\frac{2}{p}-1\right)  \right\} \setminus h(B).
$$
If we now define for some  $ \varphi = const. \in L^2_{cus}(\Gamma_{M_{\bP}}\backslash X_{\bP})$
(note that $L^2_{cus}(\Gamma_{M_{\bP}}\backslash X_{\bP}) = 
	L^2(\Gamma_{M_{\bP}}\backslash X_{\bP})$ as $\Gamma_{M_{\bP}}\backslash X_{\bP}$
	is compact), $\varphi\neq 0,$
$$
 F: \Omega \to L^p(M),\, z\mapsto E(\bP| \varphi, h(z)\rho_{\bP})
$$
the map $F_f: \Omega\to \C, z\mapsto \int_M F(z)(x)f(x) dx$ is analytic as a composition of analytic mappings for all $f\in L^{p'}(M), \frac{1}{p} + \frac{1}{p'} = 1$. Note that the integral is always finite as the 
Eisenstein series $F(z)$ are contained in $L^p(M)$. Furthermore, it follows from
Corollary \ref{eigenvalues minimal} and Theorem \ref{theorem eigenvalues} that each 
$F(z)$ is an eigenfunction of $(\DMp-c)$
for the eigenvalue $z$ and the proof of part (a) is complete.\\

If $p\geq 2$, the point spectrum of $\DMp$, and hence of $(\DMp - c)$, is a discrete subset of $\R$.
On the other hand, the intersection of the point spectrum of the generator of a chaotic semigroup 
with the imaginary axis is always infinite, cf.  \cite{MR1855839}.
\end{proof}
%
%
\subsubsection{Periods of $\boldsymbol{L^p}$ heat semigroups}
\begin{definition}
If $T(t)$ denotes a strongly continuous semigroup on a Banach space ${\cal B}$,
we call any $t>0$ such that there is some $f\neq 0$ with $T(t)f=f$ and $T(s)f\neq f$ for all $0<s<t$
a {\em period} of the semigroup $T(t)$.   
\end{definition}
In \cite{Munoz:wj} it is shown that the periods of a strongly continuous semigroup can be determined if
the eigenvalues of the generator on the imaginary axis are known:
\begin{lemma}\label{characterization periods}
Let $T(t)$ denote a strongly continuous semigroup on a Banach space ${\cal B}$ with generator $A$. Then $t>0$ is a period of $T(t)$ if and only if there exist $\alpha_k\in i\R \cap \sigma_{pt}(A), k=1,\ldots,l$ ($l=\infty$ is allowed) 
such that $t$ is the smallest positive number with
$$
t\alpha_k \in 2\pi i \Z
$$ 
for all $k$.
\end{lemma}

We use this result in order to describe the set of periods of the semigroups
$e^{-t(\DMp -c)}: L^p(M) \to L^p(M)$  on a non-compact locally symmetric space 
$M=\Gamma\backslash X$
 with arithmetic fundamental group $\Gamma$.

Let $\bP$ denote a minimal rational parabolic subgroup.
As the boundary of the parabolic region
$$
 {\cal P}_{M,p}(\rho_{\bP}) = \left\{ ||\rho_{\bP}||^2 - z^2 : 
 			z\in\C, |\re z|\leq ||\rho_{\bP}||\cdot |\frac{2}{p}-1|\right\}\subset \C,
$$
is parametrized by the curve
$$
 s\mapsto \left(\frac{2||\rho_{\bP}||}{p}+ is \right)\left(2||\rho_{\bP}|| - \frac{2||\rho_{\bP}||}{p} - is\right), 
  \qquad s\in\R,
$$
the intersection of ${\cal P}_{M,p}(\rho_{\bP}) - c$ with the imaginary axis consists of the points
in $i[-r,r]$ where $r= r(c,p) = 2||\rho_{\bP}||\left(1- \frac{2}{p}\right)\sqrt{c-c_p}$ if $c>c_p$.
\begin{proposition}
Let $M=\Gamma\backslash X$ denote a non-compact locally symmetric space with
arithmetic $\Gamma$.
\begin{itemize}
 \item[\textup{(a)}] If $p\in (1,2)$ and $c>c_p$ all but finitely many points in 
 		             $[2\pi/r(c,p),\infty)$ are periods of the semigroup
		             $e^{-t(\DMp -c)}: L^p(M) \to L^p(M)$.
  \item[\textup{(b)}] If $p\geq 2$ the semigroup $e^{-t(\DMp -c)}: L^p(M) \to L^p(M)$
  			     has no periods.
\end{itemize}
\end{proposition}
\begin{proof}
This follows immediately from Lemma \ref{characterization periods} and Theorem \ref{theorem eigenvalues}.
\end{proof}
\paragraph*{Acknowledgement} 
We want to thank the referee for many valuable suggestions which led to a great improvement of this paper.


\bibliographystyle{amsplain}
\bibliography{dissertation,hypercyclic,symmetricSpaces}

\end{document}